\documentstyle[amstex]{amsart}

 \newtheorem{conj}{Conjecture}[section]

 \theoremstyle{definition}
 
 \theoremstyle{remark}
 
 \numberwithin{equation}{section}

 \newcommand{\Real}{\mathbb{R}}

\begin{document}

\title[Convergence of Dirichlet series]
{The conditional convergence of the Dirichlet series of an $L$-function}

\author[M. O. Rubinstein]
{M.\ O.\ Rubinstein\\ \\
Pure Mathematics \\University of Waterloo\\200 University Ave W\\Waterloo, ON, Canada\\N2L 3G1}

\subjclass{Primary $L$-functions}



\begin{abstract}
The Dirichlet divisor problem is used as a model to give a conjecture concerning
the conditional convergence of the Dirichlet series of an $L$-function.
\end{abstract}

\maketitle

\section{Introduction}

Let
\begin{equation}
    L(s) = \sum_1^\infty \frac{b(n)}{n^s}
\end{equation}
be a Dirichlet series and let $\Re{s}=\sigma$.
A classical summation by parts gives
\begin{equation}
    \label{eq:1}
    \sum_{n \leq X} \frac{b(n)}{n^s}
    = \frac{1}{X^s}
    \sum_{n \leq X} b(n)
    + s \int_1^X \sum_{n \leq x} b(n) \frac{dx}{x^{s+1}}.
\end{equation}
Say that
\begin{equation}
    \label{eq:2}
    \sum_{n \leq X} b(n) = O(X^{\sigma_0})
\end{equation}
for some $\sigma_0 \in \Real$.
Then, for $\sigma> \sigma_0$, letting  $X \to \infty$, (\ref{eq:1}) converges
and becomes
\begin{equation}
    \label{eq:3}
    \sum_1^\infty \frac{b(n)}{n^s}
    = s \int_1^\infty \sum_{n \leq x} b(n) \frac{dx}{x^{s+1}}.
\end{equation}
Subtracting (\ref{eq:1}) from (\ref{eq:3}) and using (\ref{eq:2}) gives
a rate of convergence:
\begin{eqnarray}
    \label{eq:4}
    \sum_{n > X} \frac{b(n)}{n^s}
    &=& - \frac{1}{X^s}
    \sum_{n \leq X} b(n)
    + s \int_X^\infty \sum_{n \leq x} b(n) \frac{dx}{x^{s+1}} \notag \\
    &=& O_s(X^{\sigma_0-\sigma}).
\end{eqnarray}
It is therefore natural to ask, for a given $L(s)$, how small can we take $\sigma_0$,
i.e. in what half-plane does the Dirichlet series for $L(s)$ converge.

In the case of the Riemann zeta function or Dirichlet $L$-functions, the answer 
is immediate.
For the Riemann zeta function, $b(n)=1$, and so $\sigma_0 = 1$ is the best possible.

However, a well known trick allows one to take $\sigma_0=0$ and evaluate $\zeta(s)$ in 
the half plane $\sigma>0$ by writing
\begin{eqnarray}
    \label{eq:5}
    \zeta(s)\left(1-\frac{2}{2^{s}}\right) = 1 -\frac{1}{2^s}+\frac{1}{3^s}-\frac{1}{4^s}+\ldots
\end{eqnarray}
so that $b(n)=(-1)^{n-1}$ and $\sum_{n \leq X} b(n) = O(1)$. 

For a non-trivial Dirichlet character, $\chi(n)$, of modulus $q$,
we can take $\sigma_0=0$ for the Dirichlet series of $L(s,\chi)$
because $\sum_{n \leq X} \chi(n) = O_q(1)$.

For higher degree $L$-functions, however, the problem of obtaining 
a bound for the truncated sum of the Dirichlet coefficients is very difficult,
and the best known bounds seem to be far from the truth.
For example, let
\begin{eqnarray}
    \label{eq:6}
     f(z) = \sum_1^\infty a(n) q^n, \quad q=\exp(2\pi i z)
\end{eqnarray}
be a cusp form of weight $l$ and level $N$, and
\begin{eqnarray}
    \label{eq:7}
    L_f(s) =
    \sum_1^\infty \frac{a(n)}{n^{(l-1)/2}} \frac{1}{n^s}
\end{eqnarray}
be it's corresponding $L$-function. Each $a(n)$ is normalized by $n^{(l-1)/2}$
so that the critical line is $\Re{s}=1/2$.
Hecke \cite{Hecke} proved that
\begin{eqnarray}
    \label{eq:8}
    \sum_{n\leq X} \frac{a(n)}{n^{(l-1)/2}} = O_f(X^{1/2})
\end{eqnarray}
giving $\sigma_0=1/2$. This was improved to $\sigma_0=11/24+\epsilon$ by Walfisz \cite{Walfisz},
and when combined with the Ramanujan conjecture proved by Deligne \cite{Deligne} 
one can get $\sigma_0=1/3+\epsilon$.

However, this seems to be far from the truth.
To see what might be a reasonable value for $\sigma_0$ we consider
$\zeta(s)^k$, $k$ a positive integer,
which is, in some sense, the simplest degree $k$ $L$-function.
However, this is not a typical $L$-function in that its Dirichlet
coefficients are all positive and no cancellation occurs when
they are summed. This differs from the behaviour of entire $L$-function where
one expects $\sum_{n \leq X} b(n)$ to cancel. Once one removes the
contribution from the order $k$ pole of $\zeta(s)^k$,
we conjecture that the $k$-divisor problem provides a good model for entire $L$-functions
of degree $k$. 

Let $d_k(n)$ be the Dirichlet coefficients of
\begin{eqnarray}
    \label{eq:zeta k}
    \zeta(s)^k = \sum_1^\infty \frac{d_k(n)}{n^s},
\end{eqnarray}
and define
\begin{eqnarray}
    \label{eq:Delta}
    D_k(X) = \sum_{n \leq X} d_k(n).
\end{eqnarray}
The Dirichlet coefficient $d_k(n)$ is equal to the number of ways of writing $n$ as
a product of $k$ factors.

Assume now that, in (\ref{eq:1}), $b(n)=O(n^\epsilon)$.
Perron's formula \cite{Murty}[pg 67] states, 
that:
\begin{eqnarray}
    \label{eq:9}
    \sum_{n<X} b(n) = \frac{1}{2\pi i} \int_{c-iT}^{c+iT} L(s) \frac{X^s}{s} ds
    +O\left(\frac{X^{c+\epsilon}}{T}\right),
\end{eqnarray}
where $c>1$.
In the case of $L(s) = \zeta(s)^k$, one proceeds by shifting the line integral to the left
and estimating the integral along the four sides of the resulting rectangle. This gives
\begin{eqnarray}
    \label{eq:main and remainder}
    D_k(X) =  X P_k(\log{X}) + \Delta_k(X)
\end{eqnarray}
with $P_k$ being a polynomial of degree $k-1$
coming from the residue of the $k$-th order pole at $s=1$,
and $\Delta_k(X)$ denoting the remainder term.

The $k$ divisor problem states that the true order
of magnitude for $\Delta_k$ is:
\begin{eqnarray}
    \label{eq:divisor problem}
    \Delta_k(X) = O\left(X^{(k-1)/2k+\epsilon}\right).
\end{eqnarray}
When $k=2$, the traditional Dirichlet divisor problem is
\begin{eqnarray}
    \label{eq:k=2}
    D_2(X) = X\log{X} +(2\gamma-1)X+\Delta_2(X),
\end{eqnarray}
with a conjectured remainder
\begin{eqnarray}
    \label{eq:dirichlet divisor problem}
    \Delta_2(X) = O\left(X^{1/4+\epsilon}\right).
\end{eqnarray}

The estimate (\ref{eq:divisor problem}) for the remainder term $\Delta_k(X)$
is based on expected cancellation in Voronoi-type formulas for $\Delta_k(X)$
(such as (12.4.4) and (12.4.6) described in \cite{Titch}), and also on
estimates for the mean square of $\Delta_k$. For example, it is known for
$k=2,3$ \cite{Cramer} \cite{Tong} and conjectured for $k\geq 4$, that
\begin{eqnarray}
    \label{eq:mean square}
    \frac{1}{X} \int_0^X \Delta_k^2(y) dy \sim c_k X^{(k-1)/k}
\end{eqnarray}
where $c_k>0$ is constant. For $k=4$, Heath-Brown obtained a slightly weaker
upper bound, $O(X^{3/4+\epsilon})$ rather than the asymptotic \cite{HB}.
This asymptotic is known to be equivalent to
the Lindel\"of hypothesis.
For a discussion on the Dirichlet divisor problem, see Titchmarsh \cite{Titch}[Chapters XII,XIII].

By analogy with the $k$ divisor problem, it seems reasonable to conjecture:
\begin{conj}
Let $L(s)= \sum b(n)/n^s$ be an entire $L$-function
of degree $k$, normalized so the critical line is through $\Re{s}=1/2$. Then
    \begin{eqnarray}
        \label{eq:conj}
        \sum_{n \leq X} b(n) = O\left(X^{(k-1)/2k+\epsilon}\right).
    \end{eqnarray}
More generally, let $L(s)$ be meromorphic with its only pole being at $s=1$ of order $r$.
Then
    \begin{eqnarray}
        \label{eq:conjb}
        \sum_{n \leq X} b(n) = X P_L(\log{X})+ O\left(X^{(k-1)/2k+\epsilon}\right)
    \end{eqnarray}
where $P_L$ is a polynomial of degree $r-1$.
\end{conj}

This conjecture has been stated in various specific cases
\cite{Ivic} \cite{Fomenko},
but the author could not find a reference that mentions this conjecture for a general
$L$-function.

Notice that $(k-1)/2k < 1/2$, hence, if this conjecture is true, then
the Dirichlet series of an entire $L$-function can be used to evaluate it at
any $s \in \mathbb{C}$, by summing the terms for $\sigma \geq 1/2$, and by
applying the functional equation for $\sigma <1/2$. 

To mention just two examples, the Dirichlet series of a cusp form, $L_f(s)$ in (\ref{eq:7}),
is expected
to converge for $\sigma>1/4$, and the Dirichlet series of 
its symmetric square $L$-function, which has degree $3$, should converge 
for $\sigma>1/3$.

The rate of convergence, however, makes this
not very practical for smaller values of $\sigma$. For instance, 
when $k=3$, summing $10^6$ terms of the Dirichlet series gives, using
(\ref{eq:4}), less
than four decimal places accuracy at $s=1$,
while at $s=1+100i$ one expects about two digits accuracy.
To get $16$ digits at $s=1$ would be impossible from a practical point of view, 
requiring roughly $10^{24}$ terms of the Dirichlet series. 

Nonetheless, it is interesting to know, in principle and also as a way to
double check more sophisticated algorithms for computing $L$-functions, that
the Dirichlet series of an entire $L$-function does converge up to and slightly
beyond its critical line.

\subsection*{Acknowledgment}
 
I would like to thank Aleksandar Ivi\'c for comments on this paper.
This work is supported by an NSERC Discovery grant and also
by NSF FRG grant DMS-0757627.

\bibliographystyle{amsplain}


\end{document}